\documentclass[11pt]{amsart}
\usepackage{amssymb}
\usepackage{pstricks}
\psset{unit=1pt} \psset{arrowsize=4pt 1} \psset{linewidth=1pt}
\newrgbcolor{mygreen}{.2 .7 .2}
\newrgbcolor{myred}{.7 .3 .2}
\newgray{mygray}{.90}
\countdef\x=23 \countdef\y=24 \countdef\z=25 \countdef\t=26

\def\hdom(#1,#2)#3{
\x=#1 \y=#2 \multiply\x by 16 \multiply\y by 16 \z=\x \t=\y
\advance\z by 32 \advance\t by 16
\psline(\x,\y)(\x,\t)(\z,\t)(\z,\y)(\x,\y) \advance\x by 16
\advance\y by 8 \rput(\x,\y){{\bf #3}}}

\def\vdom(#1,#2)#3{
\x=#1 \y=#2 \multiply\x by 16 \multiply\y by 16 \z=\x \t=\y
\advance\z by 16 \advance\t by 32
\psline(\x,\y)(\x,\t)(\z,\t)(\z,\y)(\x,\y) \advance\x by 8
\advance\y by 16 \rput(\x,\y){{\bf #3}}}

\newtheorem{thm}{Theorem}
\newtheorem{proposition}[thm]{Proposition}
\newtheorem{lem}[thm]{Lemma}

\def\spin{{\rm spin}}
\def\weight{{\rm weight}}
\def\g{{\mathcal G}}
\def\sign{{\rm sign}}
\def\nv{{\rm nv}}
\def\bv{{\rm bv}}
\def\r{{\rm r}}

\author[Thomas Lam]{Thomas Lam}
\address{Department of Mathematics, Harvard University, Cambridge MA, 02138, USA}
\email{tfylam@math.harvard.edu}

\urladdr{http://www.math.harvard.edu/~tfylam}

\begin{document}

\title{On Sj\"{o}strand's skew sign-imbalance identity}
\begin{abstract}
Recently, Sj\"{o}strand gave an identity for the sign-imbalance of
skew shapes.  We give a quick proof of this using the skew domino
Cauchy identity and some sign analysis for skew shapes.
\end{abstract}
\maketitle

\section{The Theorem}
Let $T$ be a standard Young tableau with skew shape $\lambda/\mu$.
We will use the English notation for our tableaux, so that
partitions are top-left justified.  The {\it reading word} $r(T)$ of
$T$ is obtained by reading each row from left to right, starting
with the top row.  The {\it sign} $\sign(T)$ of $T$ is the sign of
$\r(T)$ as a permutation.  The {\it sign-imbalance}
$I_{\lambda/\mu}$ of $\lambda/\mu$ is given by
$$
I_{\lambda/\mu} = \sum_{T} \sign(T),
$$
where the summation is over all tableaux $T$ with shape
$\lambda/\mu$.

For a partition $\lambda$, let $v(\lambda) = \sum_i \lambda_{2i}$
denote the sum of the even parts. denote the Generalizing an earlier
conjecture of Stanley~\cite{Sta}, Sj\"{o}strand~\cite{S} proved the
following identity.

\begin{thm}[\cite{S}]
\label{thm:S} Let $\alpha$ be a partition and let $n \in {\mathbb
N}$ be even. Then
$$
\sum_{\lambda/\alpha \vdash n} (-1)^{v(\lambda)}
I_{\lambda/\alpha}^2 = \sum_{\alpha/\mu \vdash n} (-1)^{v(\mu)}
I_{\alpha/\mu}^2.
 $$
\end{thm}

The aim of this note is to give a quick derivation of
Theorem~\ref{thm:S} using the techniques developed in~\cite{Lam} and
the {\it skew domino Cauchy identity}.  Let $\g_{\lambda/\mu}(X;q) =
\sum_D q^{\spin(D)}x^{\weight(D)}$ be the spin-weight generating
function of domino tableaux with shape $\lambda/\mu$; see for
example~\cite{Lam}.  Here we will use the convention that $\spin(D)$
is equal to {\it half} the number of vertical dominoes in $D$.
Though not stated explicitly, the following identity is a
straightforward generalization of the ``domino Cauchy identity''
proved in any of~\cite{Lam,Lam1,vL}.

\begin{thm}
\label{thm:skewcauchy} Let $\alpha, \beta$ be two fixed partitions.
Then
$$
\sum_{\lambda} \g_{\lambda/\alpha}(X;q) \g_{\lambda/\beta}(Y;q) =
\prod_{i,j}\frac{1}{(1-x_iy_j)(1-qx_iy_j)} \sum_{\mu}
\g_{\alpha/\mu}(X;q) \g_{\beta/\mu}(Y;q).
$$
\end{thm}

\section{The Proof}
Let $D$ be a standard domino tableau with shape $\lambda/\mu$.  The
sign $\sign(D)$ is equal to $\sign(T)$ where $T$ is the standard
Young tableau obtained from $D$, also with shape $\lambda/\mu$, by
replacing the domino labeled $i$ by the numbers $2i-1$ and $2i$. The
following result follows from a sign-reversing
involution~\cite{Lam,S,Sta}.

\begin{lem}
\label{lem:first} If $\lambda/\mu$ has an even number of squares,
then its sign-imbalance is given by
$$
I_{\lambda/\mu} = \sum_D \sign(D),
$$
where the summation is over all standard domino tableau with shape
$\lambda/\mu$.

\end{lem}

Let $\delta \in D$ be a vertical domino occupying squares in rows
$i-1$ and $i$. We call $\delta$ {\it nice} if the number of squares
contained in $\lambda/\mu$ to the left of $D$ in row $i$ is odd.  In
other words, if $\delta$ lies in column $j$ then it is nice if and
only if $j - \mu_i$ is even.  Let $\nv(D)$ denote the number of nice
(and thus vertical) dominoes in $D$.  Let $\bv(D)$ denote the number
of non-nice vertical dominoes in $D$.  Then we have $\spin(D) =
\frac{\nv(D) + \bv(D)}{2}$.

\begin{lem}
\label{lem:nv} Let $D$ be a domino tableau of shape $\lambda/\mu$.
Then $\sign(D) = (-1)^{\nv(D)}$.
\end{lem}
\begin{proof}
This follows from the same argument as in the proof of
\cite[Proposition 21]{Lam}.
\end{proof}

\begin{figure}[ht]
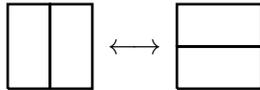

\pspicture(20,0)(70,25)  \vdom(0,0){} \vdom(1,0){}
\rput(48,15){$\longleftrightarrow$}
 \hdom(4,0){} \hdom(4,1){}
\endpspicture
\caption{The ``local move'' which connects all domino tilings.}
\label{fig:localmove}
\end{figure}

\begin{lem}
\label{lem:nvbv} Let $D$ be a domino tableau of shape $\lambda/\mu$.
The number $\nv(D) - \bv(D)$ depends only on the shape
$\lambda/\mu$.
\end{lem}
\begin{proof}
The number $\nv(D)-\bv(D)$ only depends on the tiling of
$\lambda/\mu$ by dominoes.  It is well known (see for
example~\cite{Pak}) that every two such domino tilings can be
connected by a number of moves of the form shown in
Figure~\ref{fig:localmove}. These moves do not change the value of
$\nv(D)-\bv(D)$.
\end{proof}

For each skew shape $\lambda/\mu$, define $v'(\lambda/\mu) := \nv(D)
- \bv(D)$ where $D$ is any domino tableau with shape $\lambda/\mu$.

\begin{lem}
\label{lem:mod2} We have $v(\lambda)+v(\mu) \equiv v'(\lambda/\mu)
\mod 2$.
\end{lem}
\begin{proof}
This is straight forward to prove by induction on the size of
$\lambda$, starting with $\lambda = \mu$ and adding dominoes.
\end{proof}

\begin{proposition}
\label{prop:main} Let $\lambda/\mu$ have an even number of squares.
Then
$$
I_{\lambda/\mu}^2 = (-1)^{v(\lambda/\mu)} \left(\sum_{D}
(-1)^{\spin(D)} \right)^2,
$$
where the summation is over all domino tableaux of shape
$\lambda/\mu$.
\end{proposition}
\begin{proof}
We have
\begin{align*}
I_{\lambda/\mu}^2 & = \left( \sum_D \sign(D)\right)^2 & \mbox{by
Lemma~\ref{lem:first}} \\&= \left( \sum_D (-1)^{\nv(D)} \right)^2 &
\mbox{by
Lemma~\ref{lem:nv},} \\
&= \left( \sum_D (-1)^{\spin(D) + v'(\lambda/\mu)/2}
\right)^2 & \mbox{using $\spin(D) = \frac{\nv(D) + \bv(D)}{2}$,} \\
&= (-1)^{v(\lambda/\mu)} \left(\sum_{D} (-1)^{\spin(D)} \right)^2 &
\mbox{by Lemma~\ref{lem:mod2}}.
\end{align*}
\end{proof}

\begin{proof}[Proof of Theorem~\ref{thm:S}]
In Theorem~\ref{thm:skewcauchy}, let $\alpha = \beta$ and $q = -1$.
 Then take the coefficient of $x_1x_2\cdots x_n y_1 y_2 \cdots y_n$ on
both sides.  Using Proposition~\ref{prop:main}, this gives exactly
Theorem~\ref{thm:S}.
\end{proof}

\end{document}